# A note on the WKB solutions of difference equations

M.I. Ayzatsky[1]
National Science Center Kharkov Institute of Physics and Technology (NSC KIPT),
610108, Kharkov, Ukraine

Comparison of approximate solutions that were obtained by using different asymptotic methods of solutions of difference equations with the exact solution is presented. Results show that for the studied equation the method of transformation of the difference equation into the Riccati equation gives more correct solutions than the method of direct asymptotic expansion.

## 1 Introduction

It is common knowledge (see, for example, [1,2,3]) that a second order difference equation

$$y_{k+2} + f_{1,k} y_{k+1} + f_{0,k} y_k + f_k = 0 \qquad (1)$$

has approximate solutions if the sequences $f_{0,k}$ and $f_{1,k}$ vary sufficiently slowly with $k$. As in the theory of asymptotic solutions[2] of differential equations, there are two approaches for obtaining approximate solutions of difference equations. In the frame of the first approach, we are finding of approximate solution straightforwardly of the equation (1) [3,4,5]. The second one implies transformation of the equation (1) into the Riccati equation and finding an approximate solution of this equation [6,7]. It can be shown that for differential equations these two methods are equivalent ones and they give the same results. In the case of difference equations the situation is slightly different. The method of finding an approximate solution of the Riccati equation [6] differs from the one that is used for finding an approximate solution straightforwardly of the equation (1) [3,4]. Recently, a new method of transformation of the linear difference equations into a system of the first order equations was proposed [8]. It can also be used for finding an approximate solution of difference equations.

In this note we compare approximate solutions that were obtained by using these methods and compare them with the exact solution.

## 2 The discrete WKB equations

We represent the solution of the difference equation (1) as the sum of new grid functions [8]

$$y_k = y_{1,k} + y_{2,k}. \qquad (2)$$

By introducing new unknowns $y_{n,k}$ instead of the one $y_k$, we can impose an additional condition. This condition we write in the form

$$y_{k+1} = g_{1,k} y_{1,k} + g_{2,k} y_{2,k} \qquad (3)$$

where $g_{n,k}$ ($1 \leq n \leq 2$) are the arbitrary sequences.

If

$$\det\begin{pmatrix} 1 & 1 \\ g_{1,k} & g_{2,k} \end{pmatrix} \neq 0, \qquad (4)$$

then the representation (2)-(3) is unique. Indeed, from (2) and (3) we can uniquely find $y_{n,k}$ as a linear combination of $y_k$. Using (1),(2) and (3) we can write such system of equations

---
[1] M.I. Aizatskyi, N.I.Aizatsky; aizatsky@kipt.kharkov.ua
[2] So called the LG-WKB solutions

$$\begin{pmatrix} y_{1,k} \\ y_{2,k} \end{pmatrix} = T_k \begin{pmatrix} y_{1,k-1} \\ y_{2,k-1} \end{pmatrix} + F_k. \qquad (5)$$

where

$$T_k = \begin{pmatrix} -\dfrac{f_{0,k-1} + g_{1,k-1}(g_{2,k} + f_{1,k-1})}{g_{1,1,k} - g_{1,2,k}} & -\dfrac{f_{0,k-1} + g_{2,k-1}(g_{2,k} + f_{1,k-1})}{g_{1,k} - g_{2,k}} \\ \dfrac{f_{0,k-1} + g_{1,k-1}(g_{1,k} + f_{1,k-1})}{g_{1,k} - g_{2,k}} & \dfrac{f_{0,k-1} + g_{2,k-1}(g_{1,k} + f_{1,k-1})}{g_{1,k} - g_{2,k}} \end{pmatrix}, \qquad (6)$$

$$F_k = \begin{pmatrix} -\dfrac{1}{g_{1,k} - g_{2,k}} f_k \\ \dfrac{1}{g_{1,k} - g_{2,k}} f_k \end{pmatrix} = \begin{pmatrix} -\overline{f}_k \\ \overline{f}_k \end{pmatrix} \qquad (7)$$

In some cases, it is useful to work with S-matrix (see, for example, [8,9,10,11])

$$\begin{pmatrix} y_{2,k-1} \\ y_{1,k} \end{pmatrix} = S_k \begin{pmatrix} y_{1,k-1} \\ y_{2,k} \end{pmatrix} + \overline{F}_k, \qquad (8)$$

where

$$S_k = \begin{pmatrix} -\dfrac{T_{k,21}}{T_{k,22}} & \dfrac{1}{T_{k,22}} \\ \dfrac{T_{k,11} T_{k,22} - T_{k,12} T_{k,21}}{T_{k,22}} & \dfrac{T_{k,12}}{T_{k,22}} \end{pmatrix}, \qquad (9)$$

$$\overline{F}_k = \begin{pmatrix} -\dfrac{1}{T_{k,22}} \overline{f}_k \\ -\dfrac{T_{k,12} + T_{k,22}}{T_{k,22}} \overline{f}_k \end{pmatrix}. \qquad (10)$$

We would like to emphasize that the sequences $g_{1,k}$ and $g_{2,k}$ are the arbitrary ones, and we do not impose a condition that the new grid functions $y_{1,k}, y_{2,k}$ are the solutions of the equation (1).

We will consider the homogeneous difference equations. It is means that $f_k = 0$.

Now we can describe several properties of the normal system of difference equations

$$\begin{pmatrix} y_{1,k} \\ y_{2,k} \end{pmatrix} = T_k \begin{pmatrix} y_{1,k-1} \\ y_{2,k-1} \end{pmatrix}. \qquad (11)$$

From (6) it follows that we can choose the sequences $g_{1,k}$ and $g_{2,k}$ in a such way that matrix $T_k$ will be triangular or even diagonal one. It is realized by setting $T_{k,12} = 0$ and $T_{k,21} = 0$.

These conditions give the non-linear second-order rational difference equation (Riccaty type [8,12,13,14]) for the sequences $g_{1,k}$ and $g_{2,k}$

$$f_{0,k} + g_{(1,2),k}(g_{(1,2),k+1} + f_{1,k}) = 0. \qquad (12)$$

In this case the matrix $T_k$ is a diagonal one and the system (11) takes the form:

$$\begin{aligned} y_{1,k} &= g_{1,k-1} y_{1,k-1} \\ y_{2,k} &= g_{2,k-1} y_{2,k-1} \end{aligned}. \qquad (13)$$

Solutions $y_{1,k}, y_{2,k}$ are the linearly independent ones.

The characteristic equation of the difference equation (1) is
$$\rho_k^2 + f_{1,k}\rho_k + f_{0,k} = 0. \tag{14}$$

Let $g_{1,k} = \rho_k^{(1)}$, $g_{2,k} = \rho_k^{(2)}$, where $\rho_k^{(1)}$, $\rho_k^{(2)}$ are the solutions of the characteristic equation (14)

$$\rho_k^{(1,2)} = -\frac{f_{1,k}}{2} \pm \frac{1}{2}\sqrt{f_{1,k}^2 - 4f_{0,k}}. \tag{15}$$

The matrix $T_k$ takes the form
$$T_k = \begin{pmatrix} \rho_{k-1}^{(1)} \dfrac{\rho_{k-1}^{(1)} - \rho_k^{(2)}}{\rho_k^{(1)} - \rho_k^{(2)}} & \rho_{k-1}^{(2)} \dfrac{\rho_{k-1}^{(2)} - \rho_k^{(2)}}{\rho_k^{(1)} - \rho_k^{(2)}} \\ \rho_{k-1}^{(1)} \dfrac{\rho_k^{(1)} - \rho_{k-1}^{(1)}}{\rho_k^{(1)} - \rho_k^{(2)}} & \rho_{k-1}^{(2)} \dfrac{\rho_k^{(1)} - \rho_{k-1}^{(2)}}{\rho_k^{(1)} - \rho_k^{(2)}} \end{pmatrix}. \tag{16}$$

If the sequences $f_{0,k}$ and $f_{1,k}$ vary sufficiently slowly with $k$ ($f_{0,k} = f_0(\varepsilon k)$, $f_{1,k} = f_1(\varepsilon k)$, $0 \leq \varepsilon \ll 1$), then the differences $\left(\rho_k^{(1,2)} - \rho_{k-1}^{(1,2)}\right)$ are the small values and we can neglect the non-diagonal terms in the matrix $T_k$. This gives

$$y_{1,k} = \rho_{k-1}^{(1)}\left(1 - \frac{\rho_k^{(1)} - \rho_{k-1}^{(1)}}{\rho_k^{(1)} - \rho_k^{(2)}}\right)y_{1,k-1}$$
$$y_{2,k} = \rho_{k-1}^{(2)}\left(1 + \frac{\rho_k^{(2)} - \rho_{k-1}^{(2)}}{\rho_k^{(1)} - \rho_k^{(2)}}\right)y_{2,k-1} \tag{17}$$

It can be shown that these equations coincide with the equations of the discrete WKB approach (see, for example, [6,7]). Indeed, from (13) it follows that the discrete WKB equations can be obtained with using an approximate solutions of the Riccati equation (12) under assumption that $f_{0,k}$ and $f_{1,k}$ vary sufficiently slowly with $k$. The Riccati equations can be transformed with an iteration procedure into the quadratic equations

$$\left(g_{(1,2),k-1}\right)^2 + g_{(1,2),k-1}f_{1,k-1} + f_{0,k-1} + \rho_k^{(1,2)}\left(\rho_{k+1}^{(1,2)} - \rho_k^{(1,2)}\right) = 0. \tag{18}$$

It is one of the possible forms of the quadratic equation (compare with [6,7]) that can be obtained at the second iteration. Its solutions differ from the ones that were obtained in [6,7] by an amount of order $\varepsilon^2$.

The approximate solutions of this equations with error of $O(\varepsilon^2)$ are

$$g_{(1,2),k-1} \approx \rho_{k-1}^{(1,2)}\left(1 \mp \frac{\rho_k^{(1,2)} - \rho_{k-1}^{(1,2)}}{\rho_k^{(1)} - \rho_k^{(2)}}\right). \tag{19}$$

Comparison (17) and (19) shows that these two different approaches give the same result, and the equations (17) and (13) coincide.

The solutions of the equations (17) at $k > k_0$ can be write as

$$y_k^{(1,2)} = \prod_{s=k_0+1}^{k} T_{s,11,22} y_{k_0}^{(1)} = y_{k_0}^{(1,2)} \prod_{s=k_0+1}^{k} \rho_{s-1}^{(1,2)} \left(1 \mp \frac{\rho_s^{(1,2)} - \rho_{s-1}^{(1,2)}}{\rho_s^{(1)} - \rho_s^{(2)}}\right) =$$

$$= y_{k_0}^{(1,2)} \exp\left(\sum_{s=k_0+1}^{k} \ln\rho_{s-1}^{(1,2)} + \sum_{s=k_0+1}^{k} \ln\left(1 \mp \frac{\rho_s^{(1,2)} - \rho_{s-1}^{(1,2)}}{\rho_s^{(1)} - \rho_s^{(2)}}\right)\right) \approx$$

$$\approx y_{k_0}^{(1,2)} \exp\left(\sum_{s=k_0+1}^{k} \ln\rho_{s-1}^{(1,2)} \mp \sum_{s=k_0+1}^{k} \frac{\rho_s^{(1,2)} - \rho_{s-1}^{(1,2)}}{\rho_s^{(1)} - \rho_s^{(2)}}\right) =$$

$$= y_{k_0}^{(1,2)} \exp\left(\sum_{s=k_0+1}^{k} \ln\rho_{s-1}^{(1,2)} - \sum_{s=k_0+1}^{k} \frac{\sqrt{f_{1,s}^2 - 4f_{0,s}} - \sqrt{f_{1,s-1}^2 - 4f_{0,s-1}}}{2\sqrt{f_{1,s}^2 - 4f_{0,s}}} \pm \sum_{s=k_0+1}^{k} \frac{f_{1,s} - f_{1,s-1}}{2\sqrt{f_{1,s}^2 - 4f_{0,s}}}\right) \quad (20)$$

We can transform the second term in the exponent into the multiplier

$$y_k^{(1,2)} \sim \frac{1}{\left(f_{1,k}^2 - 4f_{0,k}\right)^{1/4}} \exp\left(\sum_{s=k_0+1}^{k} \ln\rho_{s-1}^{(1,2)} \pm \sum_{s=k_0+1}^{k} \frac{f_{1,s} - f_{1,s-1}}{2\sqrt{f_{1,s}^2 - 4f_{0,s}}}\right). \quad (21)$$

The comparison this formula with the one that was obtained by finding an approximate solution straightforwardly of the equation (1) [3] gives some difference. The formula (21) contains additional sum in the exponent (the second sum).

### 3 Numerical model

In the frame of the Coupling Cavity Model (CCM) electromagnetic field in each cavity of the chain of resonators are represented as the expansion with the short-circuit resonant cavity modes [15,16,17,18,19,20]

$$\vec{E}^{(k)} = \sum_q e_q^{(k)} \vec{E}_q^{(k)}(\vec{r}) , \quad (22)$$

where $q = \{0, m, n\}$ and such coupling equations for $e_{010}^{(n)}$ can be obtained [21,22,23,24]

$$Z_k e_{010}^{(k)} = \sum_{j=-\infty, j\neq n}^{\infty} e_{010}^{(j)} \alpha_{010}^{(k,j)}. \quad (23)$$

Here $e_{010}^{(k)}$ - amplitudes of $E_{010}$ modes, $Z_k = 1 - \frac{\omega^2}{\omega_{010}^{(k)2}} - \alpha_{010}^{(k,k)}$, $\omega_{010}^{(k)}$ - eigen frequencies of these modes, $\alpha_{010}^{(k,j)}$ - real coefficients that depend on both the frequency $\omega$ and geometrical sizes of all volumes. Sums in the right side can be truncated

$$Z_k^{(N)} e_{010}^{(N,k)} = \sum_{j=k-N, j\neq k}^{k+N} e_{010}^{(N,j)} \alpha_{010}^{(k,j)}. \quad (24)$$

In the case of $N = 1$, the system of coupled equations (24) is very similar to the one that can be constructed on the basis of equivalent circuits approach (see, for example [25,26,27,28]). But in the frame of the CCM the coefficients $\alpha_{0mn}^{(k,j)}$ are electrodynamically strictly defined for arbitrary $N$ and can be calculated with necessary accuracy. In the theory of RF filters the coupling matrix circuit model is used intensively (see, for example, [29] and cited there literature). The main problem is how to calculate the matrix elements.

Amplitudes of other modes ($(m,n) \neq (1,0)$) can be found by summing the relevant series

$$e_{0mn}^{(k)} = \frac{\omega_{0mn}^{(k)2}}{\omega_{0mn}^{(k)2} - \omega^2} \sum_{j=k-N}^{k+N} e_{010}^{(j)} \alpha_{0mn}^{(k,j)}. \quad (25)$$

For the chain of cylindrical resonators longitudinal component of electric field at $r = 0$ (on the system longitudinal axis) is:

$$E_z^{(k)} = \sum_{m,n} e_{0mn}^{(j)} \cos\left(\frac{\pi}{d} nz\right). \tag{26}$$

If we can ignore "long coupling" interaction, the set of coupling equations (24) takes the form

$$Z_k e_{010}^{(k)} = e_{010}^{(k-1)} \alpha_{010}^{(k,k-1)} + e_{010}^{(k+1)} \alpha_{010}^{(k,k+1)}, \tag{27}$$

where $Z_k = \left(1 - \frac{\omega^2}{\omega_{010}^{(k)2}} - \alpha_{010}^{(k,k)} - i\frac{\omega}{\omega_{010}^{(k)} Q_k}\right)$

The set of coupling equations (27) can be considered as the second-order difference equation. This difference equation, which defines the amplitudes of the basic modes $e_{010}^{(k)}$, is the main equation of the CCM. It is reasonable to note that the amplitudes of the basic modes $e_{010}^{(k)}$ are non-measured values. Indeed, we can measure the components of electric field in any point, for example, by the nonresonant perturbation method, but we cannot measure $e_{0mn}^{(k)}$ and have to use numerical methods for finding these amplitudes by using the expansion (22). This circumference creates difficulties in studding the properties of the real slow-wave waveguides, including their tuning [30,31]. The similar situation arises also in other electrodynamic models. For example, the space harmonics in homogeneous periodic waveguides are non-measured values, too.

Below we will consider the chain of cylindrical resonators that are connected via circular central openings in the walls – the disk loaded waveguides (DLW)[3]. It was shown that the DLWs, that are usually used in linacs, with disk spacing large enough ( $d \geq \lambda/3$ ) can be describe with sufficient accuracy by the difference equation (27) [32]. Appropriate values of the coupling coefficients $\alpha_{010}^{(k,k)}, \alpha_{010}^{(k,k+1)}$ at fixed frequency can be approximated by some functions of geometrical sizes. Calculations on the base of the CCM show that for the most often used in linacs DLWs such approximations can be used

$$\alpha_{010}^{(k,k)} = -\frac{u_k \overline{p}_k + u_{k+1} \overline{p}_{k+1}}{\tilde{b}_k^2 \tilde{d}_k}$$

$$\alpha_{010}^{(k,k-1)} = \frac{u_k}{\tilde{b}_k^2 \tilde{d}_k}, \tag{28}$$

$$\alpha_{010}^{(k,k+1)} = \frac{u_{k+1}}{\tilde{b}_k^2 \tilde{d}_k}$$

where $u_k = \frac{\alpha a_k^3}{b_*^2 d_*} p_k^{(c)}$, $\tilde{b}_k = \frac{b_k}{b_*}$, $\tilde{d}_k = \frac{d_k}{d_*}$, $a_k$ - the hole radius between $k-1$ and $k$ resonators, $b_k$ - the radius of $k$ cylindrical resonator, $d_k$ - the resonator length, $b_* = c\frac{\lambda_{01}}{\omega}$, $d_*$ - normalizing parameters, $\omega_{010}^{(k)} = c\frac{\lambda_{01}}{b_k}$, $J_0(\lambda_{01}) = 0$, $\alpha = \frac{2}{3\pi J_1^2(\lambda_{01})}$, $\overline{p}_k = \frac{p_k^{(s)}}{p_k^{(c)}}$.

Analysis shows that we can consider parameters $p_k^{(s)}, p_k^{(c)}$ as the functions of the geometric sizes of the diaphragms only (the opening radius $a_k$, the thickness $t_k$ of the diaphragm between $k-1$ and $k$ resonators and the radius of the rounding of the disk hole edges).

For $t_k = 0.4$ cm, $d_k = 3.0989$ cm, parameters $p_k^{(s)}, p_k^{(c)}$ can be represented[4] as

---

[3] DLW structures are the most often used in linacs and represent the chain of cavities in which the phase varies smoothly from cell to cell in such way, that an accelerated particle constantly locates in accelerating field.

[4] For simplicity, we will consider the case without of the rounding of the disk hole edges. For taking into account the rounding of the disk hole edges.

$$p_k^{(s)} = 0.0142a_k^2 - 0.1329a_k + 0.9133$$
$$p_k^{(c)} = -0.0928a_k^2 + 0.4491a_k - 0.0444 \quad . \tag{29}$$

The parameter $p_k^{(c)}$ determines the deviation of the dependence of the coupling coefficient $\alpha_{010}^{(k,k-1)}$ on $a_k$ from the law $a_k^3$, $p_k^{(s)}$ - the deviation of the dependence of the resonator frequency shift due the hole in the $k$-disk on $a_k$ from the law $a_k^3$ (see (28)).

For simplicity we will suppose that $p_k^{(s)} = p_k^{(c)} = 1$. This assumption corresponds the case when the cavities are connected through the small round openings in the thin diaphragms.

We will consider the chain of the identical resonators ($\tilde{b}_k = \frac{b_k}{b_*} = 1$, $\tilde{d}_k = \frac{d_k}{d_*} = 1$) that consist of two semi-infinite homogeneous parts (I and II) that are connected via the inhomogeneous transient cells. The variable parameter is the radius of the opening that connects two resonators. The difference equation (27) can be written as ($y_k = e_{010}^{(k)}$)

$$y_{k+1} + \left[\frac{2(1-\cos\varphi_I) - (\bar{u}_k + \bar{u}_{k+1})}{\bar{u}_{k+1}}\right] y_k + \frac{\bar{u}_k}{\bar{u}_{k+1}} y_{k-1} = 0, \tag{30}$$

where $\varphi_I$ - the phase shift between cells in the first homogeneous part,

$$\bar{u}_k = \frac{u_k}{u_I} = \frac{(1-\cos\varphi_I)}{(1-\cos\varphi_k)} \tag{31}.$$

If we set $\varphi_k = \varphi_I$ ($\bar{u}_k = 1$) the equation (30) has such solutions of the characteristic equation (14)

$$\rho_k^{(I,1,2)} = -\frac{f_{1,k}}{2} \pm \frac{1}{2}\sqrt{f_{1,k}^2 - 4f_{0,k}} = \cos\varphi_I \pm i\sin\varphi_I = \exp(\pm i\varphi_I). \tag{32}$$

Using (31), we can transform the difference equation (30) as

$$y_{k+1} + \left[-2\cos\varphi_{k+1} + \frac{\cos\varphi_{k+1} - \cos\varphi_k}{(1-\cos\varphi_k)}\right] y_k + \left[1 - \frac{\cos\varphi_{k+1} - \cos\varphi_k}{(1-\cos\varphi_k)}\right] y_{k-1} = 0. \tag{33}$$

If the differences $(\cos\varphi_{k+1} - \cos\varphi_k)$ are small, then the solutions of the characteristic equation can be estimated as

$$\rho_k^{(1,2)} \approx \exp(\pm i\varphi_{k+1}). \tag{34}$$

Let consider the task of diffraction of the incident wave on the gradual transition between two homogeneous chains of resonators. In this case in the first homogeneous part we have two waves: the incident wave ($\rho_k^{(I,1)} = \exp(i\varphi_I)$) with the amplitude equals to 1 and the reflected wave ($\rho_k^{(I,2)} = \exp(-i\varphi_I)$) with the amplitude $R$ (the reflection coefficient). In the second homogeneous part we have only one transmitted wave ($\rho_k^{(II,1)} = \exp(i\varphi_{II})$) with the amplitude $T$ (the transmission coefficient).

It is useful to use the S-matrix method for finding solution of the considered task. If we find the S-matrix of the part of the infinite chain that contains the transition cells and several homogeneous cells on both sides of the inhomogeneous part (cell numbers are $k = 1 \div N$) and the set of the additional S-matricies, we get the solution of the considered task. Indeed, from the definition of the S-matrix we have

$$y_k = y_{1,k} + y_{2,k}, \tag{35}$$

$$y_{2,1} = S_{11}^{(N)} y_{1,1} + S_{12}^{(N)} y_{2,N},$$
$$y_{1,N} = S_{21}^{(N)} y_{1,1} + S_{22}^{(N)} y_{2,N}. \tag{36}$$

Suppose $y_{1,1} = 1$, $y_{2,N} = 0$, then

$$y_{2,1} = S_{11} = R,$$
$$y_{1,N} = S_{21} = T. \quad (37)$$

For the additional $S^{(k)}$-matrices

$$y_{2,1} = S_{11}^{(k)} y_{1,1} + S_{12}^{(k)} y_{2,k}, \quad 2 \le k \le N,$$
$$y_{1,k} = S_{21}^{(k)} y_{1,1} + S_{22}^{(k)} y_{2,k}. \quad (38)$$

From these equations it follows that

$$y_{2,k} = \frac{R - S_{11}^{(k)}}{S_{12}^{(k)}},$$
$$y_{1,k} = S_{21}^{(k)} + S_{22}^{(k)} \frac{R - S_{11}^{(k)}}{S_{12}^{(k)}}, \quad (39)$$

$$y_k = y_{1,k} + y_{2,k} = S_{21}^{(k)} + \left(S_{22}^{(k)} + 1\right) \frac{R - S_{11}^{(k)}}{S_{12}^{(k)}}. \quad (40)$$

If we know the transfer matrix $T_k$ we can find the S-matrix by using the formula (9). We will use three type of the transfer matrix $T_k$. The first matrix is the strict one (16)

$$T_k^{(0)} = \begin{pmatrix} \rho_{k-1}^{(1)} \frac{\rho_{k-1}^{(1)} - \rho_k^{(2)}}{\rho_k^{(1)} - \rho_k^{(2)}} & \rho_{k-1}^{(2)} \frac{\rho_{k-1}^{(2)} - \rho_k^{(2)}}{\rho_k^{(1)} - \rho_k^{(2)}} \\ \rho_{k-1}^{(1)} \frac{\rho_k^{(1)} - \rho_{k-1}^{(1)}}{\rho_k^{(1)} - \rho_k^{(2)}} & \rho_{k-1}^{(2)} \frac{\rho_k^{(1)} - \rho_{k-1}^{(2)}}{\rho_k^{(1)} - \rho_k^{(2)}} \end{pmatrix}. \quad (41)$$

The second matrix

$$T_k^{(1)} = \begin{pmatrix} \rho_{k-1}^{(1)} \left(1 - \frac{\rho_k^{(1)} - \rho_{k-1}^{(1)}}{\rho_k^{(1)} - \rho_k^{(2)}}\right) & 0 \\ 0 & \rho_{k-1}^{(2)} \left(1 + \frac{\rho_k^{(2)} - \rho_{k-1}^{(2)}}{\rho_k^{(1)} - \rho_k^{(2)}}\right) \end{pmatrix} \quad (42)$$

gives such solution in the WKB approach

$$y_k^{(1,2)} \sim \frac{1}{\left(f_{1,k}^2 - 4f_{0,k}\right)^{1/4}} \exp\left(\sum_{s=k_0+1}^{k} \ln \rho_{s-1}^{(1,2)} \pm \sum_{s=k_0+1}^{k} \frac{f_{1,s} - f_{1,s-1}}{2\sqrt{f_{1,s}^2 - 4f_{0,s}}}\right). \quad (43)$$

The third matrix

$$T_k^{(2)} = \begin{pmatrix} \rho_{k-1}^{(1)} \left(1 - \frac{\overline{\rho}_k^{(1)} - \overline{\rho}_{k-1}^{(1)}}{\rho_k^{(1)} - \rho_k^{(2)}}\right) & 0 \\ 0 & \rho_{k-1}^{(2)} \left(1 + \frac{\overline{\rho}_k^{(2)} - \overline{\rho}_{k-1}^{(2)}}{\rho_k^{(1)} - \rho_k^{(2)}}\right) \end{pmatrix}, \quad (44)$$

where $\overline{\rho}_k^{(1,2)} = \pm \frac{1}{2}\sqrt{f_{1,k}^2 - 4f_{0,k}}$, gives such solution in the WKB approach

$$y_k^{(1,2)} \sim \frac{1}{\left(f_{1,k}^2 - 4f_{0,k}\right)^{1/4}} \exp\left(\sum_{s=k_0+1}^{k} \ln \rho_{s-1}^{(1,2)}\right). \quad (45)$$

This case corresponds of approximate solution that was obtained without using the Riccati equation [3,4].

For numerical calculation we used the linear phase transition between homogeneous parts

$$\varphi_k = \begin{cases} \varphi_I, & 1 \le k \le N_h \\ \varphi_I + \dfrac{\varphi_{II} - \varphi_I}{N - 2N_h}(k - N_h), & N_h + 1 \le k \le N - N_h \\ \varphi_{II}, & N - N_h + 1 \le k \le N \end{cases} \qquad (46)$$

where $\varphi_I = \dfrac{\pi}{3}$, $\varphi_{II} = \dfrac{2\pi}{3}$, $N_h = 100$, $N = 250$.

Calculation results of are presented in Figure 1-Figure 6.

Calculations on the basis of the matrix $T_k^{(0)}$ give the following values of the reflection coefficient and the transmission coefficient: $|R_0| = 3.05\text{E-}008$, $|T_0| = 1.7316$. Approximate models give $|R_1| = 0$, $|T_1| = 1.7395$, $|R_2| = 0$, $|T_2| = 1.7330$

We see that the amplitude distributions are nearly the same for three types of the transition matrix. We can also see that the phase distribution for the third matrix $T_k^{(2)}$ differs from the phase distributions for the matrices $T_k^{(0)}$ and $T_k^{(1)}$. This difference is not small.

Let estimate the difference between the distributions that are based on using the matrices $T_k^{(1)}$ and $T_k^{(2)}$. From (43) and (45) it follows that the difference is in the additional sum that stands in the exponent

$$\Delta P = \sum_{s=k_0+1}^{k} \frac{f_{1,s} - f_{1,s-1}}{2\sqrt{f_{1,s}^2 - 4 f_{0,s}}} . \qquad (47)$$

Using the estimation (34), we have

$$\Delta P \approx \sum_{s=k_0+1}^{k} \frac{\cos\varphi_s - \cos\varphi_{s-1}}{2i \sin\varphi_s} = \sum_{s=k_0+1}^{k} \frac{\Delta x_s}{2i\sqrt{1-x_s^2}} \approx$$
$$\approx \int_{x_1}^{x_k} \frac{dx}{2i\sqrt{1-x^2}} = \frac{1}{2i} \arcsin x \Big|_{x_1}^{x_k} = -i\frac{1}{2}(\varphi_{II} - \varphi_I) = -i\frac{\pi}{6} \qquad (48)$$

Results of the numerical calculations are in good agreement with this estimation (see Figure 6, cavity numbers >150).

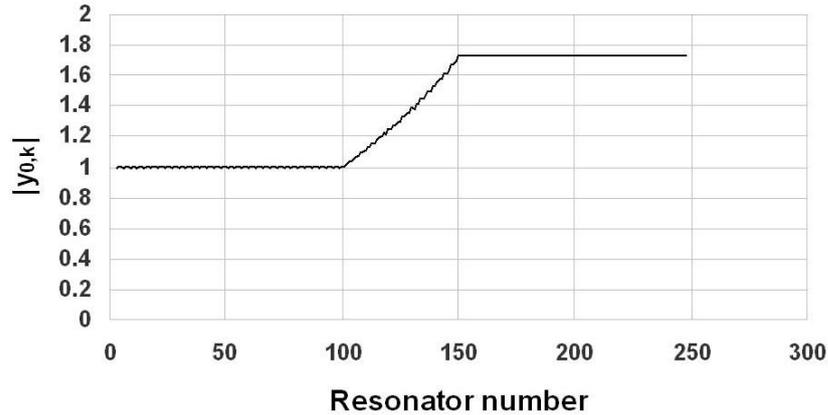

**Figure 1 The amplitude distribution of the exact solution (calculation with the matrix $T_k^{(0)}$)**

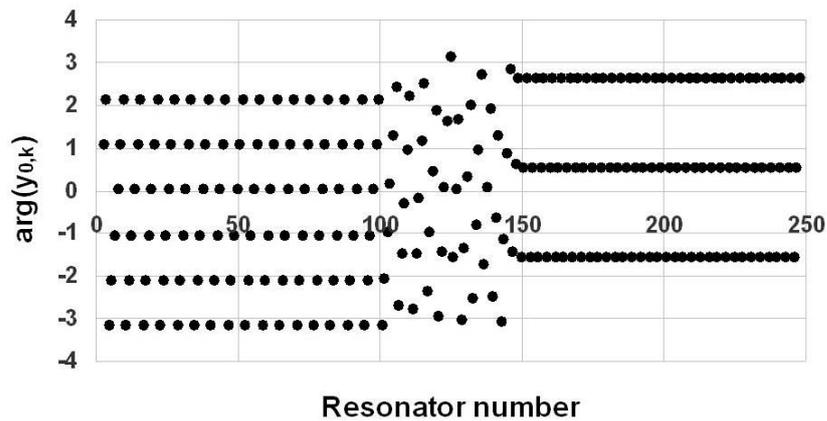

**Figure 2 The phase distribution of the exact solution**

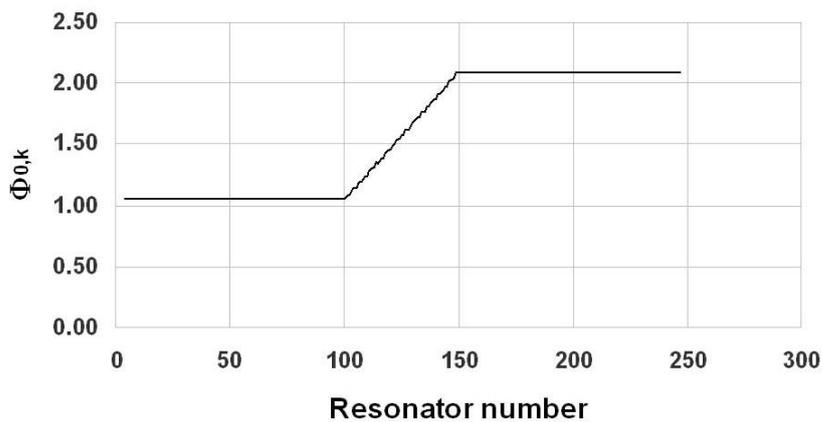

**Figure 3 The phase shift ($\Phi_{0,k} = \arg(y_{0,k+1}) - \arg(y_{0,k})$) distribution of the exact solution**

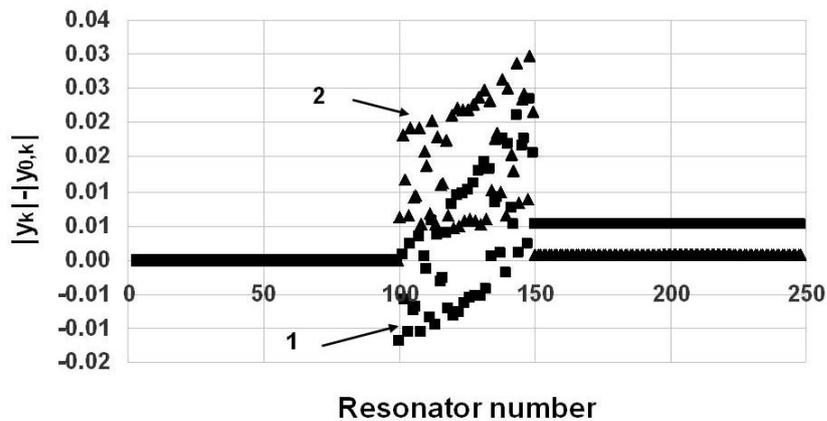

**Figure 4** Deviation of the absolute values of the approximate solutions from the exact one ($|y_k| - |y_{0,k}|$), **1**- calculation with the matrix $T_k^{(1)}$, **2** - with the matrix $T_k^{(2)}$

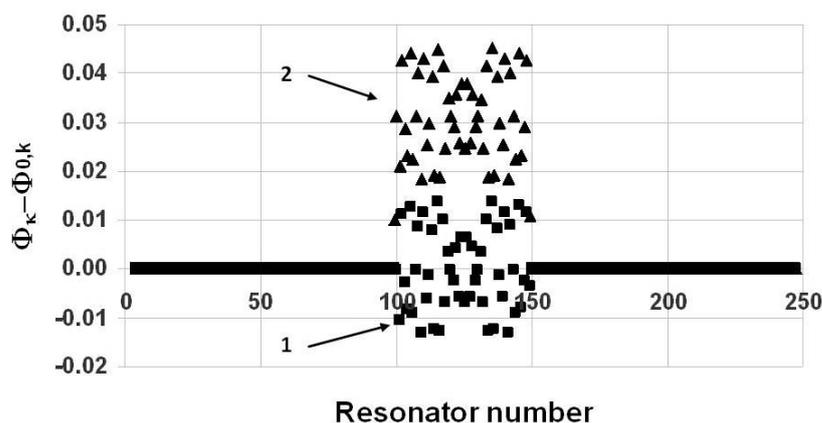

**Figure 5** Deviation of the phase shift of the approximate solutions from the exact one: 1- calculation with the matrix $T_k^{(1)}$, **2** - with the matrix $T_k^{(2)}$

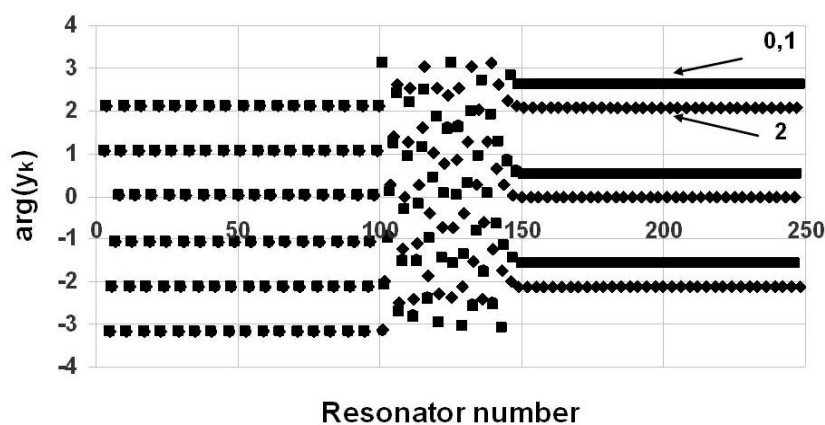

**Figure 6** The phase distribution of the different solutions: 0 - calculation with the matrix $T_k^{(0)}$, 1- with the matrix $T_k^{(1)}$, **2** - with the matrix $T_k^{(2)}$

## Conclusions

We compared approximate solutions that were obtained by using different asymptotic methods of solutions of difference equations with the exact solution. Results show that for the studied equation the method of transformation of the difference equation into the Riccati equation gives more correct solutions than the method of direct asymptotic expansion.

## Acknowledgements

The author thanks M.H.Holmes for advice and discussion.

## References

1 Spigler R., Vianello M. A survey on the Liouville-Green (WKB) approximation for linear difference equations of the second order. Advances in Difference Equations. Procceedings of the Second International Conference on Difference Equations, 1995. Gordon and Breach Science Publishers, 1997, pp.567-577

2 Garoufalidis S., Geronimo S. Asymptotics of q-difference equations. Contemporary Mathematics, V.406, Primes and Knots, AMS, 2006, pp.83-114

3 Holmes M.H. Introduction to perturbation methods Springer, 2013

4 Wilmott P. A Note on the WKB Method for Difference Equations. IMA Journal of Applied Mathematics, 1985, V. 34, pp.295-302

5 Costin O., Costin R. Rigorous WKB for finite-order linear recurrence relations with smooth coefficients. SIAM J. MATH. ANAL., 1996, V.27, N.1, pp.110-134

6 Braun P. A. WKB method for three-term recursion relations and quasienergies of an anharmonic oscillator, TMF, 1978,V.37, N.3, pp.355–370

7 Geronimo J.S. Smith D.T. WKB (Liouville-Green) Analysis od Second Order Difference Equations and Applications. Journal of Approximation Theory 1992, V.69, pp,269-301

8 M.I. Ayzatsky On the matrix form of second-order linear difference equations. http://lanl.arxiv.org/ftp/arxiv/papers//1703/1703.09608.pdf, LANL.arXiv.org e-print archives, 2017

9 Ruey-Bing (Raybeam) Hwang Periodic structures : mode-matching approach and applications in electromagnetic engineering, John Wiley & Sons Singapore Pte.Ltd., 2013

10 Akkermans E., Dunne G. V., Levy E.. Wave Propagation in One Dimension: Methods and Applications to Complex and Fractal Structures. In: Negro L.D. (ed.) Optics of aperiodic structures. Fundamentals and Device Applications. Taylor & Francis Group, 2014

11 Dobrowolski Janusz A. Microwave Network Design Using the Scattering Matrix. ARTECH HOUSE, 2010

12 Kocic V.L., Ladas G. Global behavior of nonlinear difference equations of higher order with applications. Kluwer, 1993

13 Kulenovic Mustafa R.S., Ladas G. Dynamics of second order rational difference equations. Chapman and Hall_CRC, 2001

14 Egorov A.I. Riccaty Equation. Pensoft Publishers, 2007

15 H.A. Bathe Theory of Diffraction by Small Holes. Phys. Rev., 1944, v.66, N7, p.163-182

16 V.V. Vladimirsky. ZhTF, 1947, v.17, N.11, p.1277-1282. Владимирский В.В. Связь полых электромагнитных резонаторов через малое отверстие. ЖТФ, 1947, т.17, №11, с.1277-1282

17 A.I. Akhiezer, Ya.B. Fainberg. UFN, 1951, v.44, N.3, p.321-368.Ахиезер А.И., Файнберг Я.Б. Медленные волны. УФН, 1951, т.44, №3, с.321-368

18 R.M. Bevensee. Electromagnetic Slow Wave Systems. John Wiley\&Sons, Inc.,New York-London-Sydney, 1964

19 M.A.Allen, G.S.Kino On the Theory of Strongly Coupled Cavity Chains lRE Transactions on Microwave Theory and Techniques, 1960, V.8, N.3, pp.362-372

20 R.Helm. Computer study of wave propagation, beam loading and beam blowup in the SLAC accelerator. SLAC-PUB-218, 1966

21 M.I.Ayzatsky. New Mathematical Model of an Infinite Cavity Chain. Proceedings of the EPAC96, 1996,v.3, p.2026-2028; On the Problem of the Coupled Cavity Chain Characteristic Calculations. http://xxx.lanl.gov/pdf/acc-phys/9603001.pdf, LANL.arXiv.org e-print archives, 1996

22 M.I.Ayzatskiy, K.Kramarenko. Coupling coefficients in inhomogeneous cavity chain Proceedings of the EPAC2004, 2004, pp.2759-2761

23 M.I. Ayzatskiy, V.V. Mytrochenko. Coupled cavity model based on the mode matching technique. http://lanl.arxiv.org/ftp/arxiv/papers/1505/1505.03223.pdf, LANL.arXiv.org e-print archives, 2015

24 M.I. Ayzatskiy, V.V. Mytrochenko. Coupled cavity model for disc-loaded waveguides. http://lanl.arxiv.org/ftp/arxiv/papers/1511/1511.03093.pdf, LANL.arXiv.org e-print archives, 2015

25 D.E. Nagle, E.A. Knapp, B. C. Knapp. Coupled Resonator Model for Standing Wave Accelerator Tanks. The Review of Scientific Instruments, 1967, V.38, N.11, pp.1583-1587 26 H.J.Curnow A General Equivalent Circuit for Coupled.Cavity Slow-Wave Structures. IEEE transactions on microwave theory and techniques, 1965, V.13, N. 5,pp.671-675

27 Y.Yamazaki Stability of the standing-wave accelerating strucfure studied with a coupled resonator model. Particle Accelerators, 1990, V.32, pp. 39-44

28 D.H. Whittum Introduction to Electrodynamics for Microwave Linear Accelerators. In: S.I.Kurokawa, M.Month, S.Turner (Eds) Frontiers of Accelerator Technology, World Scientific Publishing Co.Pte.Ltd., 1999

29 Miraftab Vahid, Ming Yu. Generalized Lossy Microwave Filter Coupling Matrix Synthesis and Design Using Mixed Technologies. IEEE transactions on microwave theory and techniques, 2008, V.6, N.12, pp.3016-3027

30 T.Khabiboulline, V.Puntus, M.Dohlus et al. A new tuning method for traveling wave structures. Proceedings of PAC95, pp.1666-1668

31 M.I. Ayzatskiy, V.V. Mytrochenko. Electromagnetic fields in nonuniform disk-loaded waveguides. PAST, 2016, N.3, pp.3-10; M.I. Ayzatskiy, V.V. Mytrochenko. Electromagnetic fields in nonuniform disk-loaded waveguides. http://lanl.arxiv.org/ftp/arxiv/papers/1503/1503.05006.pdf, LANL.arXiv.org e-print archives, 2015

32 M.I.Ayzatsky, E.Z.Biller. Development of Inhomogeneous Disk-Loaded Accelerating waveguides and RF-coupling. Proceedings of Linac96, p.119-121